\documentclass[11pt,a4paper]{article}
\usepackage[T1]{fontenc}
\usepackage[margin=24mm]{geometry}
\usepackage{amsmath}
\usepackage{booktabs,array}
\usepackage{listings}
\usepackage{xcolor}
\usepackage[authoryear,round]{natbib}
\usepackage{hyperref}
\hypersetup{colorlinks=true,linkcolor=blue!45!black,citecolor=blue!45!black,
urlcolor=blue!45!black,pdftitle={FunSearch and Structured SAT Search for Weak Schur Partitions},
pdfauthor={Rebecca Vicenzi, Pablo Christofferson, Caleb Hartman, Jingxuan Feng, Suleman Khan, Omar M. Elshenawy, Eric Ramos}}
\newcommand{\WS}{\mathrm{WS}}

\newcommand{\repo}{https://github.com/rebeccavicenzi/Weak_Schur_Problem_FunSearch}

\newtheorem{proposition}{Proposition}
\newenvironment{proof}{\par\noindent\textit{Proof.}\ }{\hfill\rule{1ex}{1ex}\par}
\title{FunSearch and Structured SAT Search\\for Weak Schur Partitions}
\author{Rebecca Vicenzi\textsuperscript{1} \quad
Pablo Christofferson\textsuperscript{2} \quad
Caleb Hartman\textsuperscript{3}\\
Jingxuan Feng\textsuperscript{4} \quad Suleman Khan\textsuperscript{5}\\
Omar M. Elshinawy\textsuperscript{6} \quad Eric Ramos\textsuperscript{7}\\[6pt]
{\small \textsuperscript{1}University of Bologna}\\
{\small \textsuperscript{2}Carleton College}\\
{\small \textsuperscript{3}Waynesburg University}\\
{\small \textsuperscript{4}University of Southern California}\\
{\small \textsuperscript{5}Rochester Institute of Technology}\\
{\small \textsuperscript{6}University of Münster}\\
{\small \textsuperscript{7}Stevens Institute of Technology}}
\date{}
\begin{document}
\maketitle

\begin{abstract}
The weak Schur problem asks how large an interval of consecutive integers can be partitioned into
a fixed number of classes with no class containing a monochromatic solution to $a+b=c$ with
$a\neq b$. We investigate whether FunSearch, an evolutionary program-search method guided by a
large language model, can discover effective construction strategies for this classical
combinatorial problem. In our adaptation, FunSearch evolves a \texttt{priority} function used by
a fixed greedy solver, while a separate validity check enforces the weakly sum-free condition.
This approach matches the known exact values for two, three, and four classes and reaches
$\WS(5)\geq195$. We then explore a complementary SAT-based approach. Guided by structural
information from a 207-element construction communicated by Fred Rowley
\citep{rowley2026communication}, our SAT search produces and verifies a \emph{distinct} valid
five-class partition of $\{1,\ldots,207\}$. Thus the strongest construction obtained in this
project reaches $\WS(5)\geq207$, eleven integers beyond the best published lower bound of 196.
The contrast between the two approaches also exposes a limitation of our FunSearch formulation:
its greedy solver makes local, irreversible decisions, whereas SAT can enforce global structure
across the whole partition.
\end{abstract}

\begingroup
\small
\setlength{\parskip}{0pt}
\setcounter{tocdepth}{1}
\tableofcontents
\endgroup
\newpage

\section{Introduction}
\label{sec:intro}
The weak Schur problem asks how large an interval $\{1,\ldots,N\}$ can
be while admitting a partition into a fixed number of classes, each
free of solutions to $a+b=c$ with $a\ne b$ inside a single class.
Prior work has used constructive combinatorial arguments and
computational search, including SAT-based methods
\citep{eliahou2012weak,heule2018schur,ageron2022new}.

FunSearch \citep{romeraparedes2024funsearch} uses a large language model
to evolve a short, automatically scored program that constructs a
mathematical object. We apply this idea by fixing a simple greedy
constructor and evolving the \texttt{priority} function that decides
which class to try for each new integer. A separate validity check
ensures that accepted placements preserve weak sum-freeness.

The final FunSearch program constructs a verified five-class partition
of $[1,195]$, establishing $\WS(5)\ge195$. A complementary structured
SAT search constructs a verified partition of $[1,207]$. The latter
is the strongest result obtained in this work and exceeds by eleven
the 196 lower bound reported by \citet{eliahou2012weak}.

This paper presents two computational constructions with their
algorithms and explicit, independently checkable certificates.
The structured SAT construction establishing $\WS(5)\ge207$ is the
principal result; the FunSearch construction at 195 provides the
complementary greedy approach.

\section{The Weak Schur Problem}
\label{sec:problem}
A set of positive integers $A$ is \emph{sum-free} if it contains no
solution to $a+b=c$, including equal summands $a=b$. It is
\emph{weakly sum-free} if it contains no such solution with distinct
summands. Thus the weak condition is
\begin{equation}
 a,b\in A,\quad a<b\quad\Longrightarrow\quad a+b\notin A.
 \label{eq:weak}
\end{equation}
For example, $\{1,2,4,8\}$ is weakly sum-free but not sum-free:
$1+1=2$ is allowed by the weak condition. Together with
$\{3,5,6,7\}$, it gives a weakly sum-free two-class partition of
$[1,8]$.

The weak Schur number $\WS(k)$ is the largest $N$ for which
$[1,N]=\{1,\ldots,N\}$ can be partitioned into $k$ weakly sum-free
classes. One valid partition establishes a lower bound $\WS(k)\ge N$.
An exact value additionally requires proving that no partition of
$[1,N+1]$ exists. The known small values are
\[
 \WS(1)=2,\qquad\WS(2)=8,\qquad\WS(3)=23,\qquad\WS(4)=66
\]
\citep{eliahou2012weak,heule2018schur}. For five classes,
\citet{eliahou2012weak} exhibited a partition of $[1,196]$;
subsequent template methods developed further weak Schur constructions
\citep{ageron2022new}.

Our task is to construct five lists of integers whose entries cover a
long consecutive interval exactly once, with each list weakly
sum-free. These are the class lists used by our greedy constructor
and returned by decoding the SAT output. A class label can also be
called a colour; a forbidden sum is then a monochromatic triple.

\section{FunSearch}
\label{sec:funsearch}
\textbf{FunSearch searches over programs that construct solutions.}
Here, the model proposes Python priority rules that a fixed solver
uses to construct weak Schur partitions.

FunSearch improves candidate functions through repeated generation,
execution and scoring \citep{romeraparedes2024funsearch}. A fixed
evaluator measures each construction, and promising programs become
examples for later model proposals.

The roles in our application are as follows.
\begin{description}
\item[Fixed.] The greedy solver, the validity check for a placement,
and the evaluator that scores the resulting construction.
\item[Evolved.] Only the function \texttt{priority}, which expresses
which class the solver should try first.
\item[Candidate program.] An executable Python implementation of
\texttt{priority}.
\item[Model input.] Examples of earlier candidate priority functions,
a description of the construction task, and the function interface:
its arguments and expected return value.
\item[Model output.] A new candidate priority function, written as
Python source code.
\item[Execution.] The candidate replaces \texttt{priority} in the fixed
specification and is run inside the greedy solver.
\item[Scoring.] The solver builds a partition. The fixed evaluator
assigns the candidate a numerical score from that construction.
\item[Feedback.] High-scoring candidate programs become examples in
later prompts, so subsequent proposals build on successful rules.
\end{description}

\begin{figure}[ht]
\centering
\setlength{\unitlength}{1pt}
\begin{picture}(350,195)
\put(40,165){\framebox(225,24){large language model}}
\put(144,148){\makebox(16,20){$\downarrow$}}
\put(40,127){\framebox(225,24){candidate \texttt{priority()} function}}
\put(144,110){\makebox(16,20){$\downarrow$}}
\put(40,89){\framebox(225,24){\small greedy solver builds a partition using it}}
\put(144,72){\makebox(16,20){$\downarrow$}}
\put(40,51){\framebox(225,24){partition $\longrightarrow$ numerical score}}
\put(144,34){\makebox(16,20){$\downarrow$}}
\put(40,13){\framebox(225,24){\small best-scoring candidates selected for prompts}}
\put(265,25){\rule{50pt}{0.4pt}}
\put(314.8,25){\rule{0.4pt}{152pt}}
\put(265,173){\makebox(16,8){$\leftarrow$}}
\put(279,177){\rule{36pt}{0.4pt}}
\end{picture}
\caption{The FunSearch feedback loop. The return arrow carries selected
programs back as examples for the next model prompt. The task description
and function interface remain part of that prompt.}
\label{fig:funsearch-loop}
\end{figure}
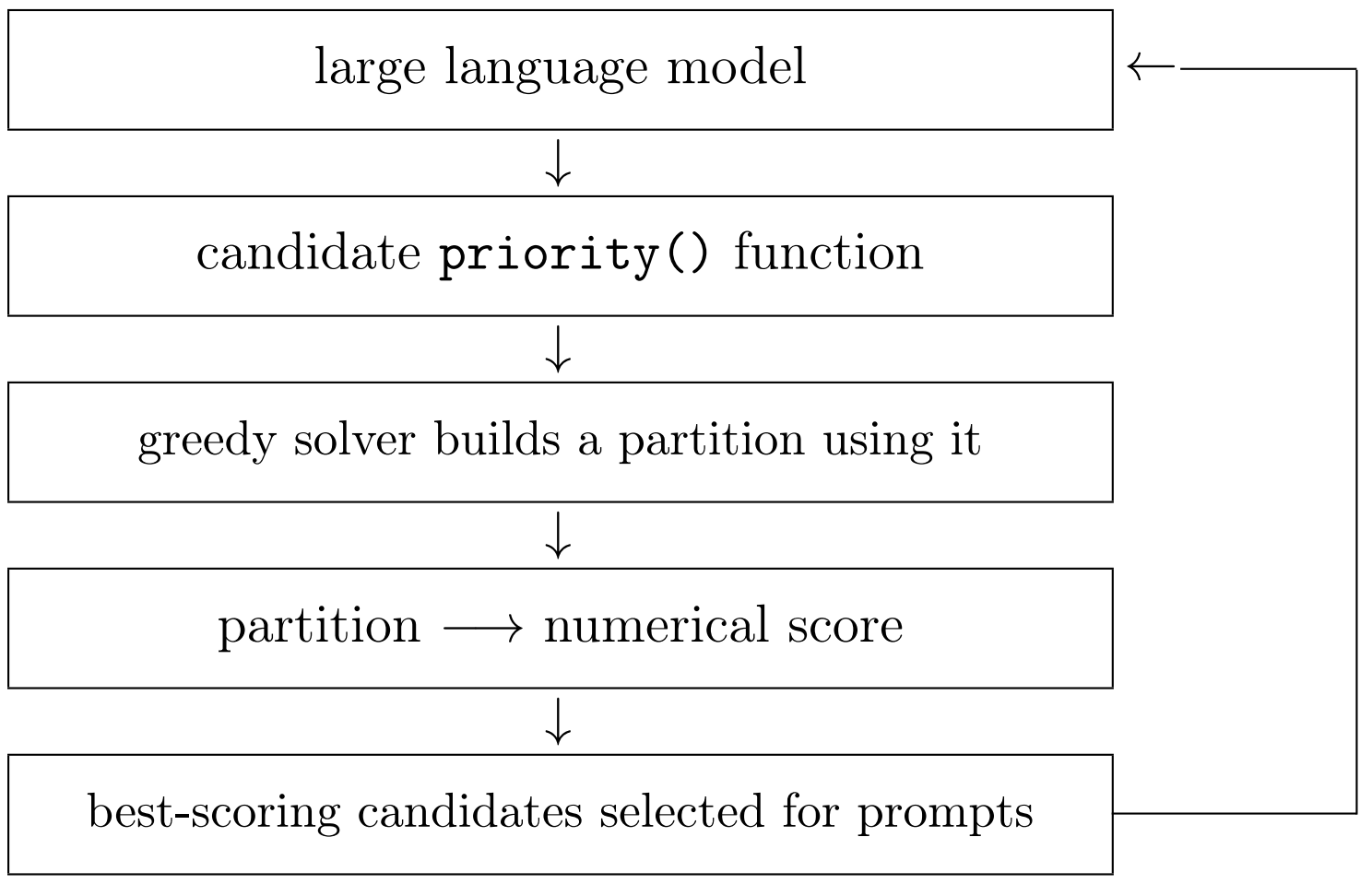

Feedback carries candidate code selected by its construction score.
Thus each rule is evaluated through a complete construction.
Section~\ref{sec:adaptation} specifies the fixed solver and the
final rule that reaches 195.

\section{Our FunSearch Adaptation and the 195 Construction}
\label{sec:adaptation}
The implementation uses \texttt{priority(i,k,partition)},
\texttt{solve(n)} and \texttt{evaluate(n)}. In these code interfaces,
$n$ is the number of classes, $k$ is the next integer to place, and
$i$ is a zero-based class index. The partial partition is a list of
$n$ class lists.

\paragraph{The evolved preference.}
\texttt{priority(i,k,partition)} returns a numerical preference for
placing $k$ in class $i$. A higher score makes that class an earlier
candidate. The function ranks choices; it does not decide whether a
choice is mathematically allowed.

\paragraph{The fixed construction.}
\texttt{solve(n)} starts with $n$ empty classes and considers
$k=1,2,3,\ldots$ in order. It computes the priorities, tries classes
from highest to lowest score and places $k$ in the first valid class.
Ties are resolved by the smallest index. If every class rejects $k$,
it returns the partition constructed so far.

The separate function \texttt{sum\_free\_addition} checks each attempted
placement. For an increasing construction, all existing elements of
a class $S$ are smaller than $k$. The only possible new forbidden
sum therefore has $k$ as its largest value: the placement is valid
exactly when no distinct $a,b\in S$ satisfy $a+b=k$.
The repository checker also tests $a+k=b$, which cannot occur on
this increasing trajectory.

For example, if a class is $\{1,2,4\}$, it cannot receive 5 because
$1+4=5$. It is rejected even if the priority ranks it first. The
solver then tries the next class. This is the essential division of
responsibility: the priority chooses preferences, while the checker
enforces weak sum-freeness.

\begin{lstlisting}[float=htbp,caption={Pseudocode for solve(n), showing the construction used for the reported result. Internal guards are omitted.},label={lst:greedy}]
start with n empty classes
for k = 1, 2, 3, ...:
    rank classes by priority(i, k, partition)
    try them from highest score to lowest, ties by index
    if a class accepts k under the validity check:
        place k in the first such class
    else:
        stop and return the current partition
\end{lstlisting}

Since every accepted insertion preserves validity and integers are
processed in order, the final read-only priority builds a partition
of an unbroken interval. No earlier placement is revised.

\paragraph{The fixed score.}
\texttt{evaluate(n)} counts the entries returned by \texttt{solve(n)}.
For five classes, only the baseline count 124 is mapped to zero;
the final construction receives score 195.

\paragraph{The final learned priority.}
The final rule looks at how far a class currently extends, values
near $k/2$, additive relationships involving $k$, values near $k$,
class size and the spread of its elements. It combines these signals
into one numerical preference used to rank the candidate classes;
the separate validity check still controls each placement.
Listing~\ref{lst:priority} reproduces the signature and computational
body from the saved repository file. NumPy is supplied as \texttt{np}
by the specification.

\begin{lstlisting}[float=htbp,caption={Verbatim excerpt of fourth\_run/WS5\_best.py: the signature and computational body of priority. Only the function's documentation string is omitted.},label={lst:priority}]
def priority(i: int, k: int, partition: list[list]) -> float:
  set_max = max(partition[i]) if partition[i] else 0
  arr = np.array(partition[i])
  window_count = np.sum(np.abs(arr - k/2) <= 2.2)
  conflict_penalty = np.sum(np.isin(arr[:, None] + arr, k)) * (k ** 0.48) / 2.9
  near_k = np.sum(arr >= k - 9) * np.exp(0.015 * (k - arr.min())) if arr.size > 0 else 0
  size_penalty = len(partition[i]) ** 0.33 if len(partition[i]) > 20 else 0
  diversity_bonus = np.std(arr) * 0.17 if arr.size > 2 else 0
  growth_penalty = -0.078 if set_max > 0.84 * k else 0
  return float(set_max) + 2.1 * window_count - 1.95 * conflict_penalty + 0.85 * near_k - size_penalty + diversity_bonus + growth_penalty
\end{lstlisting}

The pair-sum count includes equal summands, so an allowed doubling
can affect the preference. The coefficients are evolved; their
individual effects have not been validated by ablation experiments.

\begin{samepage}
Running the saved rule with \texttt{solve(5)} produces the following
classes, where $[u,v]$ denotes every integer from $u$ to $v$:
\begin{align*}
C_1={}&\{1,2,4,8,11,22,25,53,63,68,139,149,154,177,182,195\},\\
C_2={}&\{3,19,21,23\}\cup[5,7]\cup[50,52]\cup[64,66]\\
 &\quad\cup[136,138]\cup[150,152]\cup[179,181]\cup[192,194],\\
C_3={}&\{9,10,20\}\cup[12,18]\cup[54,62]\cup[140,148]\cup[183,191],\\
C_4={}&\{24,153,178\}\cup[26,49]\cup[155,176],\\
C_5={}&\{67\}\cup[69,135].
\end{align*}
\end{samepage}
Their sizes are $16,25,37,49,68$. Independent verification confirms
that each integer in $[1,195]$ occurs exactly once. Checking all
9,409 pairs $a<b$ with $a+b\le195$ finds no monochromatic triple.
The construction establishes
\[
 \boxed{\WS(5)\ge195.}
\]
Re-executing the final program reproduces this partition. Other saved
priorities construct intervals of lengths 8, 23 and 66 for two, three
and four classes, matching the known exact values.

Every class blocks 196: witness pairs are $(1,195)$, $(3,193)$,
$(9,187)$, $(24,172)$ and $(67,129)$ in $C_1,\ldots,C_5$,
respectively. This trajectory therefore stops. Extending the displayed
partition requires revising earlier choices, which the fixed greedy
solver cannot do. This motivates a global SAT formulation.

\section{SAT Formulation}
\label{sec:sat}
The greedy constructor commits permanently to each placement. SAT
instead searches for a complete assignment satisfying all constraints,
with the freedom to revise choices when conflicts arise.
A Boolean satisfiability problem consists of true-or-false variables
and clauses that must all hold. A clause is a disjunction: at least
one of its alternatives must be true.

For a target interval $[1,N]$, introduce $X_{v,c}$ for each integer
$v$ and colour $c\in\{1,\ldots,5\}$. It is true exactly when $v$
belongs to class $c$. The conjunction of the following clauses,
written in conjunctive normal form (CNF), describes the search.

\paragraph{Exactly one class per integer.}
Every integer must receive a colour, and two colours cannot both
be selected for the same integer. We express these requirements by
\begin{align}
 &\bigvee_{c=1}^{5}X_{v,c}, \label{eq:alo}\\
 &\neg X_{v,c}\lor\neg X_{v,d}\qquad(1\le c<d\le5).
 \label{eq:amo}
\end{align}
Here $\lor$ means ``or'' and $\neg$ means negation. Together the
clauses choose exactly one class for each $v$.

\paragraph{Weak Schur constraints.}
For each pair $a<b$ with $a+b\le N$, its three values cannot all
receive the same colour. For each colour $c$, add
\begin{equation}
 \neg X_{a,c}\lor\neg X_{b,c}\lor\neg X_{a+b,c}.
 \label{eq:triple}
\end{equation}
For instance, $2+3=5$ forbids giving 2, 3 and 5 the same colour.
Equal summands contribute no clause. Exactly-one and weak Schur
clauses together are equivalent to a valid five-class partition:
a partition satisfies them, and any satisfying assignment decodes
into such a partition.

\paragraph{Colour symmetry.}
The names of the five classes are arbitrary, so many SAT assignments
describe the same partition and differ only by a permutation of
colour labels. Ordering colours by their first appearance removes
this redundant symmetry: colour $c$ may appear only after colour
$c-1$ has already appeared. For $c=2,\ldots,5$ and $1\le v\le N$, add
\begin{equation}
 \neg X_{v,c}\lor\bigvee_{u=1}^{v-1}X_{u,c-1}.
 \label{eq:order}
\end{equation}
At $v=1$ the disjunction over earlier positions is empty, so only
colour 1 is allowed. Any base partition can be renamed by first
appearance, making this a convention on colour names.
We use this convention to implement the increasing-colour information
in Rowley's note \citep{rowley2026communication}, which refers to the
symmetry-breaking methods of \citet{heule2018schur}.

\paragraph{Structural restrictions.}
To focus the search at $N=207$, we adopt two restrictions from
Rowley's note \citep[p.~3 and Fig.~2]{rowley2026communication}.
These restrictions give partial colour information while leaving
the remaining assignments to the solver.

First, colours repeat with period four within $[100,120]$.
Integer $i$ and $i+4$ must have the same colour whenever both lie
in that interval. For $100\le i\le116$ and $1\le c\le5$, we impose
\begin{equation}
 \neg X_{i,c}\lor X_{i+4,c},\qquad
 \neg X_{i+4,c}\lor X_{i,c}.
 \label{eq:period}
\end{equation}
Second, fourteen colours are fixed to the values in
Table~\ref{tab:prefix}, taken from the note's Figure~2. Each becomes
a unit clause $X_{v,f(v)}$, requiring that colour to be selected.

\begin{table}[!ht]
\centering
\setlength{\tabcolsep}{4.5pt}
\begin{tabular}{r*{14}{c}}
\toprule
$v$ &1&2&3&4&5&6&7&8&9&10&11&12&14&16\\
$f(v)$ &1&2&1&3&1&2&1&3&1&2&1&4&2&3\\
\bottomrule
\end{tabular}
\caption{The fourteen prescribed colours from Rowley's structural
information.}
\label{tab:prefix}
\end{table}

The structural restrictions select a smaller family of partitions,
while retaining every base requirement. A satisfying assignment
therefore still proves a weak Schur lower bound. Failure within this
restricted family would not rule out other weak Schur partitions.

\section{The 207 Construction}
\label{sec:207}
\textbf{Structured SAT search produces a verified five-class partition
of $[1,207]$.} This is the strongest construction obtained in this
work and gives
\[
 \boxed{\WS(5)\ge207.}
\]
It exceeds by eleven the 196 lower bound reported by
\citet{eliahou2012weak}. Table~\ref{tab:main-results} summarizes the
results for the class counts relevant here.

\begin{table}[!ht]
\centering
\renewcommand{\arraystretch}{1.18}
\setlength{\tabcolsep}{5pt}
\begin{tabular}{ccc>{\raggedright\arraybackslash}p{64mm}}
\toprule
Classes $k$ & FunSearch & SAT & Published lower bound\\
\midrule
2 & 8 & -- & 8\\
3 & 23 & -- & 23\\
4 & 66 & -- & 66\\
\midrule
\textbf{5} & $\ge$\textbf{195}
& $\ge$\textbf{207}
& $\ge$\textbf{196} \citep{eliahou2012weak}\\
\bottomrule
\end{tabular}
\caption{Main computational results. FunSearch matches the known values
for two, three and four classes and reaches $\WS(5)\ge195$.
Structured SAT establishes $\WS(5)\ge207$, the strongest construction
reported here and eleven above Eliahou et~al.'s published lower bound
of 196. Dashes indicate no SAT result reported here.}
\label{tab:main-results}
\end{table}

\paragraph{The final instance and solver.}
We set $N=207$ and include exactly the constraints of
Section~\ref{sec:sat}: one class per integer, no distinct-summand
monochromatic sum, first-appearance colour ordering, period four
within $[100,120]$, and the fourteen prescribed colours. The latter
two restrictions are the partial structural information taken from
Rowley's note. Table~\ref{tab:cnf} gives the resulting clause counts.

\begin{table}[!ht]
\centering
\begin{tabular}{lr}
\toprule
Clause family & Count\\
\midrule
At least one colour &207\\
At most one colour &2,070\\
Weak Schur constraints &53,045\\
First-appearance ordering &828\\
Period four &170\\
Fixed colours &14\\
\midrule
Total &56,334\\
\bottomrule
\end{tabular}
\caption{The final SAT instance: 1,035 variables and 56,334 clauses.}
\label{tab:cnf}
\end{table}

There are 10,609 pairs $a<b$ with $a+b\le207$, each yielding five
weak Schur clauses. The generator numbers $X_{v,c}$ by $5(v-1)+c$
and writes the complete formula in DIMACS CNF format, including
the fourteen fixed colours as unit clauses.

SLIME5 \citep{riveros2020slime}, run on this file with default settings,
returns a satisfying assignment. In the recorded experiment, it found
the assignment in approximately 107 seconds. Reading each integer's
selected colour gives the five classes.

\paragraph{The certificate.}
The complete construction is supplied in the
\href{\repo/blob/f21786ee926d6535fe01767ad74392bb1c0dd985/experiments/ws5_207_sat/partition_207.json}
{machine-readable 207 certificate}. Its five lists allow the result
to be checked without repeating the search. Table~\ref{tab:results}
compares their sizes with those of the FunSearch construction.

\begin{table}[!ht]
\centering
\begin{tabular}{lrrrrrr}
\toprule
Method &$|C_1|$&$|C_2|$&$|C_3|$&$|C_4|$&$|C_5|$&$N$\\
\midrule
FunSearch &16&25&37&49&68&195\\
Structured SAT &39&41&31&46&50&207\\
\bottomrule
\end{tabular}
\caption{The independently verified five-class constructions.}
\label{tab:results}
\end{table}

\paragraph{Independent verification.}
The certificate can be checked independently by confirming that the
five classes cover $[1,207]$ exactly, with no repeated or missing
integers, and contain no forbidden weak Schur triple. Verification
examines all 10,609 pairs $a<b$ with $a+b\le207$ and confirms that
$a$, $b$ and $a+b$ never share a colour. Thus the five lists form a
partition of the required interval into weakly sum-free classes.

\begin{proposition}
The certificate partitions $[1,207]$ into five weakly
sum-free classes. Hence $\WS(5)\ge207$.
\end{proposition}
\begin{proof}
Exact interval coverage establishes a partition of $[1,207]$.
The exhaustive pair check establishes~\eqref{eq:weak} in every
class. These are the defining requirements of a weak Schur
partition, so this certificate proves the lower bound.
\end{proof}

The check is independent of the SAT solver and generator. The
\href{\repo/blob/f21786ee926d6535fe01767ad74392bb1c0dd985/experiments/ws5_207_sat/verify_partition.py}
{public verifier} equivalently checks distinct pairs within each
class and tests whether their sum is absent from that class.

The SAT assignment yields the same partition as the certificate and
satisfies all 56,334 clauses. The mathematical lower bound follows
directly from exact interval coverage and weak sum-freeness.

Rowley's source also gives a partition at 207. Our returned
certificate differs from its Figure~3 even after renaming colours:
the greatest agreement over all $5!$ colour permutations is
191 of the 207 positions.

\section{Discussion and Conclusion}
\label{sec:discussion}
FunSearch searches for a heuristic program used inside a fixed greedy
construction. The learned priority ranks candidate classes, while a
separate validity check enforces weak sum-freeness. This method
constructs a verified partition establishing $\WS(5)\ge195$, but
cannot revise earlier choices when the next insertion is blocked.
The SAT formulation performs global constraint search, allowing
assignments to be reconsidered while enforcing the full set of
requirements, and reaches $\WS(5)\ge207$.

The difference illustrates the importance of search architecture.
The gap from 195 to 207 concerns the particular greedy constructor
and learned rule used here; it does not prove a general limitation
of FunSearch. A natural future direction is a hybrid approach in
which learned heuristics guide an exact solver that can revise
assignments and enforce global constraints.

Both results are explicit, independently verified lower bounds,
rather than a determination of the exact value of $\WS(5)$.
The principal result of this paper is the structured SAT construction
establishing $\WS(5)\ge207$.

\paragraph{Code and data availability.}
The commit-stable
\href{\repo/tree/8af31072ad0cafc5946844c25eafc937be211add}{repository snapshot}
contains the
\href{\repo/blob/3424fcb034620fe81c89bcf407258a011933536e/fourth_run/WS5_best.py}{saved priority}
and
\href{\repo/blob/3424fcb034620fe81c89bcf407258a011933536e/fourth_run/WS_spec4.py}{specification},
plus the linked 207 certificate and verifier. The supplement supplies
both certificates, the CNF generator and instance, signed solver
output, independent checks and software metadata. Verification uses
only Python; re-executing the final 195 program also requires NumPy.

The final constructions can be checked from these artifacts. The
original winning FunSearch run's exact model, seed and search
configuration are not completely preserved, so the whole evolutionary
discovery process cannot be reconstructed from them.

\paragraph{Acknowledgements.}
We thank Adam Sheffer and the organizers of the 2026 Polymath Jr.
Program, where the majority of this work was completed.
We also thank Fred Rowley for providing the structural material
used in the SAT formulation.

\end{document}